\documentstyle[11pt,a4wide]{article}
\font \tenmsb=msbm10 scaled \magstep 1
\font \sevenmsb=msbm7 scaled \magstep 1
\font \fivemsb=msbm5 scaled \magstep 1
\newfam \msbfam
\textfont \msbfam=\tenmsb
\scriptfont \msbfam=\sevenmsb
\scriptscriptfont \msbfam=\fivemsb
\def \Bbb#1{\fam \msbfam \relax#1}

\font \teneufm=eufm10 scaled \magstep 1
\font \seveneufm=eufm7 scaled \magstep 1
\font \fiveeufm=eufm5 scaled \magstep 1
\newfam \eufmfam
\textfont \eufmfam=\teneufm
\scriptfont \eufmfam=\seveneufm
\scriptscriptfont \eufmfam=\fiveeufm
\def \frak#1{{\fam \eufmfam \relax#1}}

\title{\bf ON FUNCTION THEORY IN QUANTUM DISC: q-DIFFERENTIAL EQUATIONS
AND FOURIER TRANSFORM}

\author{\sl D. Shklyarov \and \sl S. Sinel'shchikov
\and \sl L. Vaksman \thanks{Partially supported by the grant
INTAS-94-4720}}

\date{\tt Institute for Low Temperature Physics \& Engineering\\
National Academy of Sciences of Ukraine}

\makeatletter
\@addtoreset{equation}{section}
\makeatother

\newtheorem{theorem}{Theorem}[section]
\newtheorem{lemma}[theorem]{Lemma}
\newtheorem{proposition}[theorem]{Proposition}
\newtheorem{corollary}[theorem]{Corollary}

\begin{document}
\maketitle

\bigskip

\section{Green function and Poisson equation}

 It was shown in \cite{SSV2} that the Laplace-Beltrami operator $\Box:L^2(d
\nu)_q \to L^2(d \nu)_q$ has a bounded inverse. Hence, for any function $f
\in L^2(d \nu)_q$, there exists a unique solution $u \in L^2(d \nu)_q$ of
Poisson equation $\Box u=f$.

\medskip

\begin{proposition}\label{deli} $\Box^{-1}f_0=-(1-q^2)\displaystyle \sum
\limits_{m=1}^\infty \frac{q^{-2}-1}{q^{-2m}-1}(1-zz^*)^m$.
\end{proposition}

\smallskip

 {\bf Proof.} It was shown in \cite[section 5]{SSV2} that the `radial part'
$\Box^{(0)}:L^2(d \nu)_q \to L^2(d \nu)_q$ of the Laplace-Beltrami operator
$\Box$ is given by $\Box^{(0)}=Dx(q^{-1}x-1)D$, with $x=(1-zz^*)^{-1}$.
Hence, $\Box^{-1}f_0=\psi(x)$,
\begin{equation}\label{dp}
\left \{{x(q^{-1}x-1)D \psi(x)=q^{-1}-q \atop \sum \limits_{j=0}^\infty
|\psi(q^{-2j})|^2 \cdot q^{-2j}<\infty}\right..
\end{equation}

 Thus, for all $x \in q^{-2{\Bbb Z}_+}$ one has
$$(q^{-2}x-1)(\psi(q^{-2}x)-\psi(x))=(q^{-1}-q)^2,$$
\begin{equation}\label{psi}
\psi(x)=\psi(q^{-2}x)-(q^{-2}-1)^2 \frac{q^4x^{-1}}{1-q^2x^{-1}}.
\end{equation}

 Now use (\ref{dp}) and (\ref{psi}) to get
$$\psi(x)=-(q^{-2}-1)^2 q^2\sum_{j=1}^\infty
\frac{q^{2j}x^{-1}}{1-q^{2j}x^{-1}}=-(q^{-2}-1)^2 q^2\sum_{j=1}^\infty
\sum_{m=1}^\infty q^{2jm}x^{-m}=$$
$$=-(q^{-2}-1)^2 q^2\sum_{m=1}^\infty
\frac{q^{2m}}{1-q^{2m}}(1-zz^*)^m.\eqno \Box$$

\medskip

 Consider the integral operator $I_m:D(U)_q \to D(U)_q'$ with the kernel
$G_m \in D(U \times U)_q'$ given by
\begin{equation}\label{G_m}
G_m=\{\left((1-\zeta \zeta^*)(1-z^*\zeta)^{-1}\right)^m \cdot
\left((1-z^*z)(1-z \zeta^*)^{-1}\right)^m \}.
\end{equation}
 The following statement was announced in \cite[Theorem 3.5]{SSV1}
\begin{theorem}\label{1} For all $f
\in D(U)_q$ \begin{equation}\label{di} \Box^{-1}f=- \sum_{m=1}^\infty
\frac{q^{-2}-1}{q^{-2m}-1}I_mf.  \end{equation}
\end{theorem}

\smallskip
 To prove this theorem we need the following auxiliary result

\begin{lemma}\label{2} $G_m$ is an invariant of the $U_q
\frak{sl}_2$-module $D(U \times U)_q'$.
\end{lemma}

\smallskip

 {\bf Proof of lemma.} The following invariants were introduced in
\cite{SSV3}:  $$k_{22}^{-m}k_{11}^{-m}=$$ $$=q^{2m}\left \{(1-\zeta
\zeta^*)^m \cdot\sum_{j=0}^\infty
\frac{(q^{2m};q^2)_j}{(q^2;q^2)_j}(q^{2(-m+1)}z^*\zeta)^j\cdot
\sum_{n=0}^\infty \frac{(q^{2m};q^2)_n}{(q^2;q^2)_n}(q^{-2m}z \zeta^*)^m
(1-z^*z)^m \right \}.$$
By a virtue of the q-binomial theorem (see \cite{GR}),
$$\sum_{i=0}^\infty \frac{(q^{2m};q^2)_i}{(q^2;q^2)_i}t^i=_1 \! \!
\Phi_0(q^{2m};-;q^2,t)=(q^{2m}t;q^2)_\infty /(t;q^2)_\infty=(t;q^2)_m^{-1}.$$
Hence,
$$k_{22}^{-m}k_{11}^{-m}=q^{2m}\left \{(1-\zeta
\zeta^*)^m(q^{-2(m-1)}z^*\zeta;q^2)_m^{-1}\cdot(q^{-2m}z
\zeta^*;q^2)_m^{-1}\cdot(1-zz^*)^m \right \}.$$

 On the other hand, in ${\rm Pol}({\Bbb C})_q$ one has $(1-\zeta
\zeta^*)\zeta=q^2 \zeta(1-\zeta \zeta^*)$, and in ${\rm Pol}({\Bbb
C})_q^{\rm op}$, respectively, $z(1-z^*z)=q^2(1-z^*z)z$, whence
$$k_{22}^{-m}k_{11}^{-m}=q^{2m}\left \{((1-\zeta
\zeta^*)(1-z^*\zeta)^{-1})^m((1-z^*z)(1-z \zeta^*)^{-1})^m \right \}.$$
The invariance of $G_m$ follows from the invariance of
$k_{22}^{-m}k_{11}^{-m}$. \hfill $\Box$

\medskip

{\bf Proof of theorem \ref{1}.} In the
special case $f=f_0$ one has $I_mf_0=(1-q^2)(1-zz^*)^m$ since
$\zeta^*f_0=f_0 \zeta=0$, $\displaystyle \int \limits_{U_q}f_0d \nu=1-q^2$.
Hence in that special case (\ref{di}) follows from proposition \ref{deli}.

 By \cite[proposition 3.9]{SSV2} $f_0$ generates the $U_q
\frak{sl}_2$-module $D(U)_q$. What remains is to show that the linear
operators $\Box^{-1}$ and $-\displaystyle \sum_{m=1}^\infty
\frac{\textstyle q^{-2}-1}{\textstyle q^{-2m}-1}I_m$ are morphisms of $U_q
\frak{sl}_2$-modules. For the first operator this follows from
\cite[proposition 4.3]{SSV2} and for the second one from lemma \ref{2}.
\hfill $\Box$

\bigskip

\section{Cauchy-Green formula}

 Let $f \in D(U)_q$. This section presents a solution of the
$\overline{\partial}$-problem in $L^2(d \mu)_q$:
\begin{equation}\label{odp}
{\partial^{(r)}\over \partial z^*}u=f,\qquad u \perp{\rm
Ker}\left({\partial^{(r)}\over \partial z^*}\right).
\end{equation}

 Our aim is to prove the following statement (see \cite[proposition 4.1]
{SSV1})

\medskip

\begin{theorem}\label{ccf} \hspace{-.5em}. Let $f \in D(U)_q$.
Then \begin{enumerate}

\item There exists a unique solution $u \in L^2(d \mu)_q$ of the
$\overline{\partial}$-problem $\overline{\partial}u=f$, which is orthogonal
to the kernel of $\overline{\partial}$.

\item $u= {\textstyle 1 \over \textstyle 2 \pi i}\displaystyle \int
\limits_{U_q} d \zeta {\textstyle \partial^{(l)}\over
\textstyle \partial z}G(z,\zeta) fd \zeta^*$, with $G
\in D(U \times U)_q^\prime$ being the Green function of the Poisson
equation.

\item $f=-{\textstyle 1 \over \textstyle 2 \pi i}\displaystyle \int
\limits_{U_q}(1-z \zeta^*)^{-1}(1-q^{-2}z \zeta^*)^{-1}d \zeta f(\zeta)d
\zeta^*- {\textstyle 1 \over \textstyle 2 \pi i}\displaystyle \int
\limits_{U_q} d \zeta{\textstyle \partial^{(l)}\over \textstyle \partial
z}G(z,\zeta) \cdot{\textstyle \partial^{(r)}f \over \textstyle \partial
\zeta^*}d \zeta^*$.

\end{enumerate}
\end{theorem}

 To clarify the symmetry of this problem, pass from the partial derivative
to the differential, and from functions to differential forms.

 Consider the morphism of $U_q \frak{su}(1,1)$-modules
$\overline{\partial}:\Omega(U)_q^{(1,0)}\to \Omega(U)_q^{(1,1)}$. By a
virtue of the canonical isomorphisms of covariant $D(U)_q$-bimodules
$\Omega(U)_{-2,q}^{(0,j)}\simeq \Omega(U)_q^{(1,j)}$, $j=0,1$,
$fv_{-2}\mapsto fdz$, $f \in \Omega(U)_q^{(0,*)}$, the following scalar
products are $U_q \frak{su}(1,1)$-invariant (see \cite{SSV2}):
$$(f_1dz,f_2dz)=\int \limits_{U_q}f_2^*f_1(1-zz^*)^2d \nu,\qquad
(f_1dzdz^*,f_2dzdz^*)=\int \limits_{U_q}f_2^*f_1(1-zz^*)^4d \nu.$$

 The completions of pre-Hilbert spaces $\Omega(U)_q^{(1,0)}$,
$\Omega(U)_q^{(1,1)}$, are canonically isomorphic to the Hilbert spaces
$L^2(d \mu)_q$, $L^2((1-zz^*)^2d \mu)_q$, respectively ($i_0:fdz \mapsto f$;
\ $i_1:fdzdz^* \mapsto f$ are just those isomorphisms).

 We may reduce solving the problem (\ref{odp}) to solving the following
problem:
\begin{equation}\label{od}
\overline{\partial}u=fdzdz^*,\qquad u \perp{\rm Ker}(\overline{\partial}),
\end{equation}
where the orthogonality means that the above invariant scalar product in the
space of $(1,0)$-forms vanishes.

 To solve this problem, we need auxiliary linear operators
$\overline{\partial}^*$,
$\Box^{(1,1)}=-\overline{\partial}\cdot \overline{\partial}^*$. Turn to
studying these operators.

\medskip

\begin{lemma}\label{dlf} For all $f \in D(U)_q$, ${\textstyle
\partial^{(l)}f^*\over \textstyle \partial z^*}=\left({\textstyle
\partial^{(r)}f \over \textstyle \partial z}\right)^*$.
\end{lemma}

\smallskip

 {\bf Proof.} $dz^*\cdot{\textstyle \partial^{(l)}f^*\over \textstyle
\partial z^*}=\overline{\partial}f^*=(\partial f)^*=\left({\textstyle
\partial^{(r)}f \over \textstyle \partial z}\cdot dz \right)^*=dz^*
\cdot \left({\textstyle \partial^{(r)}f \over \textstyle \partial
z}\right)^*$.\hfill $\Box$

\medskip

\begin{lemma}\label{ods} For all $f_1,f_2 \in D(U)_q$,
$\left(\overline{\partial}(f_1dz),f_2dzdz^*\right)=
\left(f_1dz,q^{2}{\textstyle \partial^{(r)} \over \textstyle \partial z}(f_2
\cdot(1-zz^*)^2)dz \right)$.
\end{lemma}

\smallskip

 {\bf Proof.} An application of lemma \ref{dlf} and the q-analogue of
Green's formula (see appendix in \cite{SSV1}) allows one to get for all
$f_1,f_2 \in D(U)_q$:
$$\left(\overline{\partial}(f_1dz),f_2dzdz^*\right)=-q^2\int
\limits_{U_q}f_2^*{\textstyle \partial^{(r)}f_1 \over \textstyle \partial
z^*}(1-zz^*)^2d \mu=-q^2\int \limits_{U_q}(1-zz^*)^2f_2^*{\textstyle
\partial^{(r)}f_1 \over \textstyle \partial z^*}d \mu=$$
$$={q^2 \over 2i\pi}\int
\limits_{U_q}dz(1-zz^*)^2f_2^*\overline{\partial}f_1={-q^2 \over 2i\pi}\int
\limits_{U_q}dz \overline{\partial}((1-zz^*)^2f_2^*)f_1=q^2\int
\limits_{U_q}{\textstyle \partial^{(l)}\over \textstyle \partial
z^*}((1-zz^*)^2f_2^*)f_1d \mu=$$
$$=q^2\int \limits_{U_q}\left({\textstyle \partial^{(r)}\over \textstyle
\partial z}\left(f_2(1-zz^*)^2 \right)\right)^*f_1d
\mu=q^2\left(f_1dz,{\textstyle \partial^{(r)} \over \textstyle \partial
z}(f_2 \cdot(1-zz^*)^2)dz \right).\eqno \Box$$

\medskip

\begin{corollary} The linear operator
$$\overline{\partial}^*:\Omega(U)_q^{(1,1)}\to \Omega(U)_q^{(1,0)};\qquad
\overline{\partial}^*:fdzdz^* \mapsto q^2{\textstyle \partial^{(r)} \over
\textstyle \partial z}(f \cdot(1-zz^*)^2)dz,$$
is a morphism of $U_q \frak{sl}_2$-modules.
\end{corollary}

\medskip

\begin{corollary}\label{b11} The linear operator
$\Box^{(1,1)}:\Omega(U)_q^{(1,1)}\to \Omega(U)_q^{(1,1)}$ given by
$\Box^{(1,1)}:fdzdz^* \mapsto q^4{\textstyle \partial^{(r)} \over \textstyle
\partial z^*}{\textstyle \partial^{(r)} \over \textstyle \partial
z}(f(1-zz^*)^2)dzdz^*$, $f \in D(U)_q$, is an endomorphism of $U_q
\frak{sl}_2$-modules.
\end{corollary}

\medskip

 The relation $\Box^{(1,1)}=-\overline{\partial}\cdot \overline{\partial}^*$
allows one to get a solution of the $\overline{\partial}$-problem in the
form $u=-\overline{\partial}^*\omega$, with $\omega$ being a solution of the
Poisson equation $\Box^{(1,1)}\omega=fdzdz^*$.

 Find a solution of the latter equation.

\medskip

\begin{lemma}\label{gen} The elements $\{z^m \}_{m>0}$, $z^*z$,
$\{z^{*m}\}_{m>0}$, generate the $U_q \frak{sl}_2$-module ${\rm Pol}({\Bbb
C})_q$.
\end{lemma}

\smallskip

 {\bf Proof}\ reduces to reproducing the argument used while proving
\cite[theorem 3.9]{SSV2}. \hfill $\Box$

\medskip

\begin{lemma}\label{b'} The linear operator $\Box':D(U)_q'\to D(U)_q'$ given
by $\Box':f \mapsto q^4 \left({\textstyle \partial^{(r)} \over \textstyle
\partial z^*}{\textstyle \partial^{(r)} \over \textstyle \partial z}f
\right)(1-zz^*)^2$, is an endomorphism of the $U_q \frak{sl}_2$-module
$D(U)^{\prime}_q$.
\end{lemma}

\smallskip

 {\bf Proof.} Consider the isomorphism of $U_q \frak{sl}_2$-modules $i:
\Omega(U)_q^{(0,0)}\to \Omega(U)_q^{(1,1)}$ given by $i:f \mapsto f
\cdot(1-zz^*)^{-2}dzdz^*$. Obviously, $\Box'=i^{-1}\Box^{(1,1)}i$. What
remains is to refer to corollary \ref{b11}. \hfill $\Box$

\medskip

\begin{proposition}\label{bb'} $q^2\Box=\Box'$. \end{proposition}

\medskip

 Before proving this proposition, we deduce its corollaries.

\medskip

\begin{corollary} $\Box^{(1,1)}=q^2 i \Box i^{-1}$.
\end{corollary}

\medskip

\begin{corollary}\label{be} $\Box f=q^2 \left({\textstyle \partial^{(r)}
\over \textstyle \partial z^*}{\textstyle \partial^{(r)} \over \textstyle
\partial z}f \right)(1-zz^*)^2$, $f \in D(U)_q'$.
\end{corollary}

\medskip

 Since $i$ is an isometry, and $0<c_1 \le -\Box \le c_2$ (see \cite{SSV2}),
one has

\medskip

\begin{corollary}\label{ddse} $0<c_1 \le q^{-2} \overline{\partial}\cdot
\overline{\partial}^*\le c_2$.
\end{corollary}

\medskip

 Note that we have proved the boundedness of the linear map
$\overline{\partial}^*$ from the completion of $\Omega(U)_q^{(1,1)}$ to the
completion of $\Omega(U)_q^{(1,0)}$.

\medskip

 {\bf Proof of proposition \ref{bb'}.} Let $f=z^*z$. By a virtue of
\cite[lemma 5.1]{SSV2} one has $\Box(z^*z)=-q^2 \Box
x^{-1}=q^2(1-zz^*)^2=q^{-2}\Box'f$.  Thus, the relation $\Box
f=q^{-2}\Box'f$ is proved in the special case $f=z^*z$.  In the two another
special cases $f \in \{z^m \}_{m \ge 0}$, $f \in \{z^{*m}\}_{m \ge 0}$ the
above relation follows from $\Omega f=\Box f=\Box'f=0$, with $\Omega$ being
the Casimir element (see \cite{SSV2}). Hence, by virtue of lemmas \ref{gen},
\ref{b'}, the relation $\Box f=q^{-2}\Box'f$ is valid for all the
polynomials $f \in{\rm Pol}({\Bbb C})_q$. What remains is to apply the
continuity of the linear maps $\Box$, $\Box'$ in the topological vector
space $D(U)_q'$ together with the density of ${\rm Pol}({\Bbb C})_q$ in
$D(U)_q'$. \hfill $\Box$

\medskip

 The following result, together with its proof attached below, are due to
S. Klimek and A. Lesniewski \cite{KL}.

\medskip

\begin{proposition}\label{proj} Consider the orthogonal projection $P$ from
$L^2(d \mu)_q$ onto the subspace $H^2(d \mu)_q$ generated by the monomials
$\{z^m \}_{m>0}$. For all $f \in D(U)_q$ one has $Pf=\int \limits_{U_q}(1-z
\zeta^*)^{-1}(1-q^2z \zeta^*)^{-1}f(\zeta)d \mu(\zeta)$.
\end{proposition}

\smallskip

 {\bf Proof.} An application of \cite[lemma 7.1]{SSV1} and the q-binomial
theorem (see \cite{GR}) yield the following explicit expression for the
kernel of the integral operator $P$:
$$\sum_{m=0}^\infty \frac{(q^4;q^2)_m}{(q^2;q^2)_m}(z \zeta^*)^*=(q^4z
\zeta^*;q^2)_\infty \cdot(z \zeta^*;q^2)_ \infty^{-1}.\eqno \Box$$

\medskip

 {\sc Remark 2.13} Another proof of proposition \ref{proj}, which involves
no properties of q-special functions, will be presented in appendix of
\cite{SSV5}.

\medskip

{\bf Proof of theorem \ref{ccf}.}
 By corollary \ref{be}, $\Omega(U)_q^{(1,1)}$ contains a unique solution
$\omega$ of the Poisson equation $\Box^{(1,1)}\omega=fdzdz^*$. It is given
by
$$\omega=q^{-2}\left(\int \limits_{U_q}G(z,\zeta)f(\zeta)(1-\zeta
\zeta^*)^2d \nu \right)(1-zz^*)^{-2}dzdz^*,$$
with $G \in D(U \times U)_q'$, $G=-\displaystyle \sum
\limits_{m=1}^\infty \frac{q^{-2}-1}{q^{-2m}-1}G_m$, being the Green function
found in section 1.

 By lemma \ref{ods} and corollary \ref{ddse}, the $(1,0)$-form $
\left(-\displaystyle \int \limits_{U_q}\frac{\partial^{(r)}G(z,\zeta)} {
\partial z}f(\zeta)d \mu \right)dz$ is a solution of the
$\overline{\partial}$-problem (\ref{od}). Hence, the function $u=
-\displaystyle  \int \limits_{U_q} {\textstyle \partial^{(r)}G(z,\zeta)\over
\textstyle \partial z}f(\zeta)d \mu$ is a solution of the
$\overline{\partial}$-problem (\ref{odp}). Since the uniqueness of a
solution of this $\overline{\partial}$-problem is obvious, we have proved
the first two statements of theorem \ref{ccf}.

 Let $f \in D(U)_q$, and $u=-\displaystyle \int
\limits_{U_q}\frac{\partial^{(r)}G(z,\zeta)}{\partial
z}\frac{\partial^{(r)}f}{\partial \zeta^*}d \mu$ be the above solution of
the $\overline{\partial}$-problem ${\textstyle \partial^{(r)}u \over
\textstyle \partial z^*}={\textstyle \partial^{(r)}f \over \textstyle
\partial z^*}$, $u \perp{\rm Ker}\left({\textstyle \partial^{(r)}\over
\textstyle \partial z^*}\right)$. Then $u \perp H^2(d \mu)_q$, $f-u \in
H^2(d \mu)_q$, and hence $Pf=P(f-u)=f-u$. Thus, $f=u+Pf$, and by a virtue of
proposition \ref{proj},
$$f=\int \limits_{U_q}(1-z \zeta^*)^{-1}(1-q^2z \zeta^*)^{-1}f(\zeta)d
\mu(\zeta)- \int \limits_{U_q}\frac{\partial^{(r)}G(z,\zeta)}{\partial
z}\frac{\partial^{(r)}f}{\partial \zeta^*}d \mu.$$

 This relation implies the third statement of theorem
\ref{ccf}, the Green formula.

\bigskip

\section{Eigenfunctions of the operator $\Box$}

 It follows from \cite[section 5]{SSV2} that $q\Box f=\Omega f$, $f \in
D(U)_q'$, with $\Omega \in U_q \frak{sl}_2$ being the Casimir element. Our
purpose is to produce distributions $f \in D(U)_q'$ for which $\Omega
f=\lambda f$ for some $\lambda \in{\Bbb C}$ . More exactly, we shall prove
the following result (it was announced in \cite[proposition 5.1]{SSV1})

\medskip

\begin{theorem}\label{ev}\hspace{-.5em}. For all $f \in{\Bbb C}[\partial
U]_q$ the element $$u=\int \limits_{\partial U}P_{l+1}(z,e^{i \theta})f(e^{i
\theta}){d \theta \over 2 \pi}\eqno(5.3)$$ of $D(U)_q^\prime$ is an
eigenvector of $\Box$:  $$\Box u=\lambda(l)u,\quad
\lambda(l)=-\frac{(1-q^{-2l})(1-q^{2l+2})}{(1-q^2)^2}.$$
\end{theorem}
\medskip

 We start with a similar problem for the quantum cone and, as in
\cite{SSV3}, consider the spaces $F(\widetilde{\Xi})_q^{(l)}\subset
F(\widetilde{\Xi})_q$ of degree $2l$ homogeneous functions on the quantum
cone $\widetilde{\Xi}$. Impose also the notation
$F(\Xi)_q=F(\widetilde{\Xi})_q \cap D(\Xi)_q'$,
$F(\Xi)_q^{(l)}=F(\widetilde{\Xi})_q^{(l)}\cap D(\Xi)_q'$, $l \in{\Bbb C}$.

 By the construction, $F(\Xi)_q^0$ is a covariant $*$-algebra. We intend to
give its description in terms of generators and relations.

\medskip

\begin{proposition} The bilateral ideal $J \subset{\rm Pol}({\Bbb C})_q$
generated by the single element $1-zz^*=0$, is a $U_q
\frak{sl}_2$-submodule of the $U_q \frak{sl}_2$-module ${\rm Pol}({\Bbb
C})_q$.
\end{proposition}

\smallskip

 {\bf Proof} is derivable from the explicit formulae
\begin{equation}\label{act} Hz=2z,\qquad X^-z=q^{1/2},\qquad
X^+z=-q^{-1/2}z^2,
\end{equation}
\begin{equation} Hz^*=-2z,\qquad X^+z^*=q^{-1/2},\qquad
X^-z^*=-q^{1/2}z^{*2}\qquad \Box
\end{equation}

\medskip

\begin{corollary} The $*$-algebra ${\Bbb C}[\partial U]_q \simeq{\rm
Pol}({\Bbb C})_q/J$ considered in \cite{SSV1} is a covariant $*$-algebra.
\end{corollary}

\medskip

 Remind that in ${\Bbb C}[\partial U]_q$ one has $zz^*=z^*z=1$.

 An application of the relation (1.3) of \cite{SSV3} yields

\medskip

\begin{proposition} The covariant $*$-algebra ${\Bbb C}[\partial U]_q$ is
isomorphic to the covariant $*$-algebra $F(\Xi)_q^0$ as follows:
$$i_0:{\Bbb C}[\partial U]_q \to F(\Xi)_q^0,\qquad
i_0:z \mapsto qt_{11}t_{12}^{-1},\qquad i_0:z^*\mapsto t_{21}^{-1}t_{22}.$$
\end{proposition}

\medskip

 Note that the vector spaces $F(\Xi)_q^{(l)}$, $l \in{\Bbb C}$, are
covariant $F(\Xi)_q^{(0)}$-bimodules, and the vector space $F(\Xi)_q$ is a
covariant $*$-algebra. We identify the elements of ${\Bbb C}[\partial U]_q$
and their images under the embedding $i:{\Bbb C}[\partial U]_q
\hookrightarrow F(\Xi)_q$.

 Let $l \in{\Bbb C}$, $x=t_{12}t_{12}^*=-qt_{12}t_{21}$. Apply the relation
(1.3) of \cite{SSV3} to get a description of the covariant bimodule
$F(\Xi)_q^{(l)}$.

\medskip

\begin{proposition} For all $l \in{\Bbb C}$, $x^l \in F(\Xi)_q^{(l)}$ one
has
\begin{equation} zx^l=q^{2l}xz,\qquad z^*x^l=q^{-2l}x^lz^* \end{equation}
\begin{equation}\label{act1}
\cases{X^+(x^l)=q^{-3/2}\frac{\textstyle q^{-2l}-1}{\textstyle q^{-2}-1}zx^l
\cr X^-(x^l)=q^{3/2}\frac{\textstyle 1-q^{2l}}{\textstyle 1-q^2}z^*x^l \cr
H(x^l)=0}
\end{equation}
\end{proposition}

\medskip

 The covariant bimodules $F(\Xi)_q^{(l)}$ are, in particular, $U_q
\frak{sl}_2$-modules. The associated representations of $U_q \frak{sl}_2$
are called the representations of the principal series. These are
irreducible for some open dense set of $l \in{\Bbb C}$. By a virtue of
relation (5.8) of \cite{SSV2}, for those $l \in{\Bbb C}$, and hence for all
$l \in{\Bbb C}$ and all $f \in F(\Xi)_q^{(l)}$, one has
\begin{equation}\label{so}
\Omega f=\Lambda(l)f,\qquad
\Lambda(l)=\frac{(q^{-l}-q^l)(q^{-(l+1)}-q^{l+1})}{(q^{-1}-q)^2}.
\end{equation}

 Let $V^{(l)}$ be the $U_q \frak{sl}_2$-modules considered in \cite{SSV2}.
One can easily deduce from (\ref{act}), (\ref{act1}), (\ref{so}) the
following

\medskip

\begin{corollary} For all $l \in{\Bbb C}$, the linear map $i_l:V^{(l)}\to
F(\Xi)_q^{(l)}$; $i_l:X^{\pm m}e_0 \,\mapsto X^{\pm m}(x^l)$, $m \in{\Bbb
Z}_+$, are the isomorphisms of $U_q \frak{sl}_2$-modules.
\end{corollary}

\medskip

 {\bf Proof of theorem \ref{ev}.} Let us turn to a construction of
distributions $f \in D(X)_q'$ on the quantum hyperboloid, which satisfy the
equation $\Omega f=\Lambda(l)f$ for some $l \in{\Bbb C}$.

 By the results of \cite[section 6]{SSV3}, the element
$$k_{22}^lk_{11}^l \stackrel{\rm def}{=}q^{-2l}\xi^l \sum_{j=0}^\infty
\frac{(q^{-2l};q^2)_j}{(q^2;q^2)_j}(q^{2(l+1)}z^*\zeta)^j \cdot
\sum_{m=0}^\infty
\frac{(q^{-2l};q^2)_m}{(q^2;q^2)_m}(q^{2l}z \zeta^*)^m(1-z^*z)^{-l}$$
of the completion of $F(X)^{\rm op}\otimes F(\Xi)_q^{(l)}$ is an invariant.
(Here $z,z^* \in F(X)^{\rm op}$, $\zeta,\zeta^*,\xi \in F(\Xi)_q$ are the
elements given by explicit formulae in \cite[section 6]{SSV3}).

 It follows from the results of \cite[section 4]{SSV3} that the linear
functional $\eta:F(\Xi)_q^{(-1)}\to{\Bbb C}$,
$$\int_{\Xi_q}\left(\sum_{m=-\infty}^\infty a_m \zeta^m \right)\xi^{-1}d
\eta=a_0,$$
is an invariant integral. Hence, the linear integral operator
$$F(\Xi)_q^{(-l-1)}\to F(X)_q;\qquad f \mapsto
\int_{\Xi_q}\{k_{22}^lk_{11}^l \}fd \eta$$
is a morphism of $U_q \frak{sl}_2$-modules. By a virtue of (\ref{so}), for
any trigonometric polynomial $f(\zeta)\in{\Bbb C}[\partial U]_q$, the
function
$$\int_{\Xi_q}\{k_{22}^lk_{11}^l \}f \xi^{-(l+1)}d
\eta=q^{-2l}\int_{\partial U}P_{-l}(z,e^{i \theta})f(e^{i \theta}){d \theta
\over 2 \pi}$$
is an eigenfunction of the Laplace-Beltrami operator. Here $P_{-l}$ is a
q-analogue of the Poisson kernel (see \cite[section 5]{SSV1}). Now a passage
from the quantum hyperboloid $X$ to the quantum disc $U$ via the isomorphism
of $U_q \frak{sl}_2$-modules $i:D(U)_q'{{\atop \textstyle \to} \atop
{\textstyle \approx \atop}}D(X)_q'$ (see \cite{SSV3}) yields the statement
of theorem \ref{ev} \hfill $\Box$

 Denote by ${\Bbb C}[\partial U]_{q,l}$ the vector space ${\Bbb C}[\partial
U]_q$ equipped by the structure of $U_q \frak{sl}_2$-module in such a way
that the map ${\Bbb C}[\partial U]_{q,l}\to F(\Xi)_q^{(l)}$; $f(z)\mapsto
f(z)x^l$, is a morphism of $U_q \frak{sl}_2$-modules.

 An application of (\ref{act}), (\ref{act1}) gives
$$X^+f(z)=-q^{-1/2}z^2(Df)(z)+q^{-3/2}\frac{q^{-2l}-1}{q^{-2}-1}f(qz),$$
$$X^-f(z)=q^{1/2}(Df)(z)+q^{3/2}\frac{1-q^{2l}}{1-q^2}f(qz),$$
$$Hf(z)=2z{d \over dz}f(z),$$
with $D:f(z)\mapsto(f(q^{-1}z)-f(qz))/(q^{-1}z-qz)$, $f \in{\Bbb C}[\partial
U]_{q,l}$.

 Let ${\rm Re}\,l>-{1 \over 2}$. With the notation of \cite{SSV1} being
implicit, introduce a linear operator $I_l$ in ${\Bbb C}[\partial U]_{q,l}$
given by
\begin{equation}\label{i_l}
I_lf=\frac{\Gamma_{q^2}^2(l+1)}{\Gamma_{q^2}(2l+1)}\lim_{1-r^2 \in q^{2{\Bbb
Z}_+}\atop r \to 1}(1-r^2)^lb_ru,
\end{equation}
with $u=\displaystyle \int \limits_{\partial U}P_{l+1}(z,e^{i \theta})f(e^{i
\theta}){d \theta \over 2 \pi}$, $f \in{\Bbb C}[\partial U]_{q,l}$.

 Our aim now is to prove the following result (see \cite[proposition
5.3]{SSV1})

\medskip

\begin{theorem} \label{lim}\hspace{-.5em}. $I_lf=f$
\end{theorem}

\medskip
{\bf Proof.}
  The theorem will be proved if we establish the existence
of the limit in the right hand side of (\ref{i_l}) and show that $I_l$ is
the identity operator.

 Let $L \subset{\Bbb C}[\partial U]_{q,l-1}$ be the linear subspace of all
those elements $f \in{\Bbb C}[\partial U]_q$ for which the both above
statements are valid. By a virtue of \cite[lemma 5.4]{SSV1},
\begin{equation}\label{asy}
\lim_{x \in q^{-2{\Bbb Z}_+}\atop z \to +\infty}\varphi_l \left({ 1 \over
x}\right)\left/\left(\frac{\Gamma_{q^2}(2l+1)}{\Gamma_{q^2}^2(l+1)}x^l
\right)\right.=1,
\end{equation}
for ${\rm Re}\,l>-{1 \over 2}$, with $\varphi_l=\displaystyle \int
\limits_{\partial U}P_{l+1}(z.e^{i \theta})f(e^{i \theta}){\textstyle d
\theta \over \textstyle 2 \pi}$. Thus, $1 \in L$. Moreover, an application
of this lemma and the fact that the linear operator
$$j_l:{\Bbb C}[\partial U]_{q,l}\to D(U)_q';\qquad j_l:f \mapsto \int
\limits_{\partial U}P_{l+1}(z,e^{i \theta})f(e^{i \theta}){d \theta \over 2
\pi}$$
is an isomorphism of $U_q \frak{sl}_2$-modules, allows one to prove that $L$
is a submodule of the $U_q \frak{sl}_2$-module ${\Bbb C}[\partial U]_{q,l}$.
On the other hand, with $l \notin{\Bbb Z}_++{\textstyle \pi \over \textstyle
\ln(q^{-1})}{\Bbb Z}$, the $U_q \frak{sl}_2$-module ${\Bbb C}[\partial
U]_{q,l}\simeq V^{(l)}$ is simple. Hence, for $l$ as above one has $L={\Bbb
C}[\partial U]_{q,l}$, and thus the theorem is
proved.\hfill $\Box$
\medskip

\medskip

 {\sc Remark 3.8.} Let $m \in{\Bbb Z}_+$, and $\psi(x)$ be a function on
$q^{-2{\Bbb Z}_+}$ such that $z^m \cdot \psi(y^{-1})=\displaystyle \int
\limits_{\partial U}P_{l+1}(z,e^{i \theta})e^{im \theta}{\textstyle d \theta
\over \textstyle 2 \pi}$. Another way of proving the existence of the limit
in the right hand side of (\ref{i_l}) is based on producing a fundamental
system of solutions of the difference equation for $\psi(x)$. (This
difference equation is a consequence of the relation $\Omega(z^m
\psi(y^{-1}))=\Lambda(l)(z^m \psi(y^{-1}))$.) It is easy to prove the
existence of such fundamental system of solutions $\psi_1,\psi_2$ that
$$\lim_{x \to+\infty}{\psi_1(x)\over x^l}=\lim_{x \to+\infty}{\psi_2(x)\over
x^{-l-1}}=1.$$
What remains is to use the relation ${\rm Re}\,l>-{1 \over 2}$.

\bigskip

\section{Decomposing in eigenfunctions of the operator $\Box^{(0)}$}

 One can find in \cite[section 5]{SSV2} a description of the bounded linear
operator $\Box^{(0)}:f(x)\mapsto Dx(q^{-1}x-1)Df(x)$ in the Hilbert space
$L^2(d \nu)_q^{(0)}$ of such functions on $q^{-2{\Bbb Z}_+}$ that $\|f
\|=\left(\displaystyle \int \limits_1^\infty|f(x)|^2d_{q^{-2}}x
\right)^{1/2}<\infty$. That section also contains the relation (5.9) which
determines the eigenfunctions $\Phi_l(x)$ of $\Box^{(0)}$. Besides, a
unitary operator $u:L^2(d \nu)_q^{(0)}\to L^2(dm)$ that realizes a
decomposition in those eigenfunctions was constructed. Remind that $u$
could be defined by (5.10), and $dm$ is a Borel measure on a compact ${\frak
L}_0$ introduced by (5.7) in \cite{SSV2}.

 In this section, explicit formulae for eigenfunctions $\Phi_l(x)$ and the
spectral measure will be found; \cite[proposition 3.2]{SSV1} will be proved.

\medskip

\begin{proposition} $\Phi_l(x)={_3}\Phi_2 \left[{\textstyle
x,q^{-2l},q^{2(l+1)};q^2;q^2 \atop \textstyle q^2,0}\right]$.
\end{proposition}

\smallskip

 {\bf Proof.} By a virtue of \cite[corollary 5.2]{SSV1}, the distribution
$${_3}\Phi_2 \left[{\textstyle (1-zz^*)^{-1},q^{-2l},q^{2(l+1)};q^2;q^2
\atop \textstyle q^2,0}\right]\in D(U)_q'$$
is an eigenfunction of $\Box$. What remains is to apply the definition of
$\Phi_l(x)$ and the evident relation ${_3}\Phi_2 \left[{\textstyle
1,q^{-2l},q^{2(l+1)};q^2;q^2 \atop \textstyle q^2,0}\right]=1$. \hfill
$\Box$

\medskip

\begin{corollary}\label{sp} The spectrum of $\Box^{(0)}$ coincides with the
segment $\left[-\frac{\textstyle
1}{\textstyle(1-q)^2},-\frac{\textstyle 1}{\textstyle
(1+q)^2} \right]$.
\end{corollary}

\smallskip

 {\bf Proof.} It follows from \cite[section 5]{SSV1} that the continuous
spectrum of $\Box^{(0)}$ fills this segment. So we are to prove that the
discrete spectrum of $\Box^{(0)}$ is void, that is $\Phi_l \notin L^2(d
\nu)_q^{(0)}$ for ${\rm Re}\,l>-{1 \over 2}$. This can be deduced from
proposition 5.1 and lemma 5.4 of \cite{SSV1}. \hfill $\Box$

\medskip

 By corollary \ref{sp}, the carrier of $dm$ coincides with the segment
$\{l \in{\Bbb C}|\;{\rm Re}\,l=-{1 \over 2},\:0 \le{\rm Im}\,l
\le{\textstyle \pi \over \textstyle h}\}$, with $h=-2 \ln \,q$. Hence,
$$\frac{1}{(1+q)^2}\le -\Box^{(0)} \le \frac{1}{(1-q)^2}.$$
This inequality implies \cite[proposition 3.2]{SSV1}.

  We intend to obtain  an explicit formula for the kernel $G(x,\xi,l)$ of the
integral operator $(\Box^{(0)}-\lambda(l)I)^{-1}$ in $L^2(d \nu)_q^{(0)}$.
By corollary \ref{sp}, the `Green function' $G(x,\xi,l)$ is well defined and
holomorphic in $l$ for $x,\xi \in q^{-2{\Bbb Z}_+}$, ${\rm Re}\,l \ne-{1
\over 2}$.

 Remind the notation $[a]_q=(q^{-a}-q^a)/(q^{-1}-q)$, and choose the branch
of $x^l$ in the half-plane ${\rm Re}\,x>0$: $x^l=e^{\ln \,x \cdot l}$, with
$\ln \,x$ being the principal branch of the logarithm.

\medskip

\begin{lemma} With $|x|>q^{2}$, ${\rm Re}\,x>0$, the function
\begin{equation}\label{psi_l}
\psi_l(x)=x^l \cdot{_2}\Phi_1 \left({q^{-2l},q^{-2l};q^2;q^{2}x^{-1}\atop
q^{-4l}}\right)
\end{equation}
satisfies the difference equation
\begin{equation}\label{epsi}
Dx(q^{-1}x-1)D \psi_l(x)=\lambda(l)\psi_l(x).
\end{equation}
\end{lemma}

\smallskip

 {\bf Proof.} The right hand side of (\ref{psi_l}) is of the form $x^l
\displaystyle \sum_{m=0}^\infty{a_m \over x^m}$, $a_m \in{\Bbb C}$. Its
substitution into (\ref{epsi}) gives
$${a_{m+1}\over a_m}=q\frac{[l-m]_q^2}{[l-m]_q[l-1-m]_q-[l]_q[l+1]_q}=
q\frac{[l-m]_q}{[m+1]_q[m-2l]_q}=$$
$$=q^{2}\frac{(1-q^{-2l+2m})^2}{(1-q^{2(m+1})(1-q^{-4l+2m})}.$$
What remains is to use the definition of the basic hypergeometric series
${_2}\Phi_1$ (see \cite{GR}).

\medskip

 This lemma and the definition of the Green function $G(z,\xi,l)$ imply

\medskip

\begin{proposition}\hfill\\
1) For ${\rm Re}\,l>-{1 \over 2}$
\begin{equation}\label{c_1}
G(x,\xi,l)=c_1(l)\cases{\psi_l(\xi)f_l(x),&$x \le \xi$ \cr
f_l(\xi)\psi_l(x),&$x \ge \xi$}
\end{equation}
2) For ${\rm Re}\,l<-{1 \over 2}$
\begin{equation}\label{c_2}
G(x,\xi,l)=c_2(l)\cases{\psi_{-1-l}(\xi)f_l(x),&$x \le \xi$ \cr
f_l(\xi)\psi_{-1-l}(x),&$x \ge \xi$}
\end{equation}
Here $x,\xi \in q^{-2{\Bbb Z}_+}$, $c_1(l),c_2(l)\in{\Bbb C}$.
\end{proposition}

\medskip

 Find the `constants' $c_1(l),c_2(l)$.

\medskip

\begin{lemma}\label{Leib} For any two functions $u,v$ on the semi-axis
$x>0$,
$$Du(x)\cdot v(x)=D(u(x)v(qx))-qu(qx)(Dv)(qx).$$
\end{lemma}

\smallskip

 {\bf Proof.} The following q-analogue of Leibnitz formula is directly from
the definition of $D$:
$$D(u(x)v(x))=(Du)(x)\cdot v(q^{-1}x)+u(qx)(Dv)(x).$$
Replace $v(x)$ by $v(qx)$ to get
$$(Du)(x)v(x)=D(u(x)v(qx))-u(qx)D(v(qx)).$$
What remains is to apply the straightforward relation $D(v(qx))=q(Dv)(qx)$.
\hfill $\Box$

\medskip

 Let $l \in{\Bbb C}$, $x \in q^{-2{\Bbb Z}_+}$, and
$\varphi_1(x),\varphi_2(x)$ be solutions of the difference equation
\hbox{$Dx(q^{-1}x-1)D \varphi=\lambda(l)\varphi$}.

\medskip

\begin{lemma}\label{indep}
$$W(\varphi_1,\varphi_2)=
x(q^{-2}x-1)\left(\frac{\varphi_1(q^{-2}x)-\varphi_1(x)}{q^{-2}x-x}\cdot
\varphi_2(x)-\varphi_1(x)\frac{\varphi_2(q^{-2}x)-\varphi_2(x)}{q^{-2}x-x}
\right)$$
does not depend on $x \in q^{-2{\Bbb Z}_+}$.
\end{lemma}

\smallskip

 {\bf Proof.} Evidently,
$$0=(Dx(q^{-1}x-1)D \varphi_1)\varphi_2-\varphi_1(Dx(q^{-1}x-1)D
\varphi_2).$$
Hence, by a virtue of lemma \ref{Leib},
$$0=D(x(q^{-1}x-1)D \varphi_1(x)\cdot \varphi_2(qx))-D(x(q^{-1}x-1)D
\varphi_2(x)\cdot \varphi_1(qx)).$$
That is,
$$D(x(q^{-1}x-1)(D \varphi_1(x)\cdot \varphi_2(qx)-\varphi_1(qx)D
\varphi_2(x))=0.$$
Hence, $q^{-1}x(q^{-2}x-1)((D \varphi_1(q^{-1}x)\cdot
\varphi_2(x)-\varphi_1(x)(D \varphi_2)(q^{-1}x))$ is a constant. \hfill
$\Box$.

\medskip

 Let $\varphi_1(x),\varphi_2(x)$ be the eigenfunctions involved in the
formulation of the previous lemma, and set up
$$\Phi(x,\xi)=\cases{\varphi_1(x)\varphi_2(\xi),&$x \ge \xi$ \cr
\varphi_1(\xi)\varphi_2(x),&$x \le \xi$}.$$

\medskip

\begin{lemma}
$$\left.Dx(q^{-1}x-1)D
\Phi(x,\xi)\right|_{x=\xi}=\cases{\frac{W(\varphi_1,\varphi_2)}{(1-q^2)\xi}
+\lambda \Phi(\xi,\xi), &$x=\xi$ \cr \lambda \Phi(x,\xi), &$x \ne \xi$}.$$
\end{lemma}

\smallskip

 {\bf Proof.}Let $x=\xi$:
$$\left.Dx(q^{-1}x-1)D \Phi \right|_{x=\xi}=
\left.\frac{q^{-1}x(q^{-2}x-1)\frac{\Phi(q^{-2}x,\xi)-\Phi(x,\xi)}
{q^{-2}x-x}-x(x-1)\frac{\Phi(x,\xi)-\Phi(q^2x,\xi)}{x-q^2x}}
{q^{-1}x-qx}\right|_{x=\xi}=$$
$$={1
\over(q^{-1}-q)\xi}\left(q^{-1}\xi(q^{-2}\xi-1)\frac{\varphi_1(q^{-2}\xi)-
\varphi_1(\xi)}{q^{-2}\xi-\xi} \cdot \varphi_2(\xi)-q
\xi(\xi-1)\varphi_1(\xi)\frac{\varphi_2(\xi)-\varphi_2(q^2 \xi)}{\xi-q^2
\xi}+\right.$$
$$\left.+q^{-1}\xi(q^{-2}\xi-1)\varphi_1(\xi)\frac{\varphi_2(q^{-2}\xi)-
\varphi_2(\xi)}{q^{-2}\xi-\xi}-q^{-1}\xi(q^{-2}\xi-1)\varphi_1(\xi)
\frac{\varphi_2(q^{-2}\xi)-\varphi_2(\xi)}{q^{-2}\xi-\xi}\right).$$
 We did not break the equality since we have added and then subtracted from
its right hand side the same expression:
$${1 \over(q^{-1}-q)\xi}q^{-1}\xi(q^{-2}\xi-1)\varphi_1(\xi)
\frac{\varphi_2(q^{-2}\xi)-\varphi_2(\xi)}{q^{-2}\xi-\xi}.$$
Thus we get
$$\left.Dx(q^{-1}x-1)D \Phi
\right|_{x=\xi}=\frac{q^{-1}}{(q^{-1}-q)\xi}\cdot
W(\varphi_1,\varphi_2)+\lambda \cdot \Phi(\xi,\xi).$$
In the case $x \ne \xi$ the statement of the lemma is evident. \hfill
$\Box$.

\medskip

\begin{corollary}\label{c_i}\hfill \\
1) For ${\rm Re}\,l>-{1 \over 2}$, $W(\psi_l,f_l)\ne 0$,
$c_1(l)={\textstyle 1 \over \textstyle W(\psi_l,f_l)}$,\\
2) For ${\rm Re}\,l<-{1 \over 2}$, $W(\psi_{-1-l},f_l)\ne 0$,
$c_2(l)={\textstyle 1 \over \textstyle W(\psi_{-1-l},f_l)}$
\end{corollary}

\medskip

 Find $W(\psi_l,f_l)$, $W(\psi_{-1-l},f_l)$ as functions of an indeterminate
$l$. Remind the notation (see \cite[section 6]{SSV1}):
\begin{equation}\label{c(l)}
c(l)=\frac{\Gamma_{q^2}(2l+1)}{\left(\Gamma_{q^2}(l+1)\right)^2}=
\frac{(q^{2(l+1)};q^2)_\infty^2}{(q^{2(2l+1)};q^2)_\infty(q^2;q^2)_\infty}.
\end{equation}

\medskip

\begin{lemma}\label{f_l} For all $l \notin{1 \over 2}+{\Bbb Z}$,
$f_l(x)=c(l)\psi_l(x)+c(-1-l)\psi_{-1-l}(x)$.
\end{lemma}

\smallskip

 {\bf Proof.} Consider the functions $f_l,\psi_l,\psi_{-1-l}$ holomorphic in
the domain $l \notin{1 \over 2}+{\Bbb Z}$. Evidently,
$\{\psi_l,\psi_{-1-l}\}$ form the base in the vector space of solutions for
the equation \hbox{$Dx(q^{-1}x-1)D \psi=\lambda(l)\psi$} in the space of
functions on $q^{-2{\Bbb Z}_+}$.  Hence
$f_l(x)=a(l)\psi_l(x)+b(l)\psi_{-1-l}(x)$, with $a(l),b(l)$ being
holomorphic in the domain $l \notin{1 \over 2}+{\Bbb Z}$.  Let $x \in
q^{-2{\Bbb Z}_+}$ go to infinity. By a virtue of \cite[lemma 5.4]{SSV1},
$a(l)=c(l)$ for ${\rm Re}\,l>-{1 \over 2}$, and $b(l)=c(-1-l)$ for ${\rm
Re}\,l<-{1 \over 2}$. What remains is to apply the holomorphy of
$a(l),b(l),c(l),c(-1-l)$ in the domain $l \notin{1 \over 2}+{\Bbb Z}$.
\hfill $\Box$

\medskip

\begin{lemma}\label{W} $W(\psi_l,\psi_{-1-l})=[2l+1]_q$.
\end{lemma}

\smallskip

 {\bf Proof.} With $x \in q^{-2{\Bbb Z}_+}$, $x \to +\infty$, one has
$$\psi_l(x)\sim x^l,\qquad \frac{\psi_l(q^{-2}x)-\psi_l(x)}{q^{-2}x-x}\sim
\frac{q^{-2l}-1}{q^{-2}-1}x^{l-1}.$$
Hence by lemma \ref{indep},
$$W(\psi_l,\psi_{-1-l})=\lim_{x \to+\infty \atop x \in q^{-2{\Bbb
Z}_+}}x(1-q^{-2}x)\left({q^{-2l}-1 \over q^{-2}-1}-{q^{-2(-1-l)}-1 \over
q^{-2}-1}\right)x^{-2}.\eqno \Box$$

\medskip

 Lemmas \ref{f_l}, \ref{W} and corollary \ref{c_i} imply

\medskip

\begin{proposition} The constants in (\ref{c_1}) and (\ref{c_2}) are given
by
\begin{equation}\label{c12}
c_1(l)=\frac{1}{c(-1-l)[2l+1]_q},\qquad c_2(l)=-\frac{1}{c(l)[2l+1]_q},
\end{equation}
with $c(l)$ being the q-analogue of Harish-Chandra's c-function determined
by (\ref{c(l)}).
\end{proposition}

\medskip

 The conclusion is as follows. For ${\rm Re}\,l \ne -{1 \over 2}$ the
operator $\Box^{(0)}-\lambda(l)I$ in the Hilbert space $L^2(d \nu)_q^{(0)}$
has a bounded inverse operator given by
$$((\Box^{(0)}-\lambda(l)I)^{-1}\psi)(x)=\int \limits_1^\infty
G(x,\xi,l)\psi(\xi)d_{q^2}\xi,\qquad \psi \in L^2(d \nu)_q^{(0)}.$$
The Green function is given by the explicit formulae (\ref{c_1}),
(\ref{c_2}), (\ref{c12}).

 Find the spectral projections of $\Box^{(0)}$.

 The following well known result follows from the Stieltjes inversion
formula (see \cite{N}).

\medskip

\begin{proposition}\label{specpr} Let $A$ be a bounded selfadjoint operator
with simple purely continuous spectrum. For any interval $(a_1,a_2)$ on the
real axis, one has
$$E((a,b))=\lim_{\varepsilon \to+0}{1 \over 2 \pi i}\int
\limits_{a_1}^{a_2}(R_{\lambda-i \varepsilon}-R_{\lambda+i \varepsilon})d
\lambda,$$
with $R_ \lambda=($A$-\lambda I)^{-1}$.
\end{proposition}

\medskip

 {\sc Remark 4.13.} There is an extension of proposition \ref{specpr} to the
case of an arbitrary selfadjoint operator (see \cite[chapter 10, section
6]{DS}).

\medskip \stepcounter{theorem}

\begin{proposition}\label{dgdlst} Let $x,\xi \in q^{-2{\Bbb Z}_+}$, ${\rm
Re}\,l=-{1 \over 2}$. Then
\begin{equation}\label{dgdl}
\lim_{\varepsilon \to+0}(G(x,\xi,l+\varepsilon)-G(x,\xi,l-\varepsilon))=
\frac{f_l(\xi)f_l(x)}{c(l)c(-1-l)[2l+1]_q}.
\end{equation}
\end{proposition}

\smallskip

 {\bf Proof.} In the case $x \le \xi$ one has due to (\ref{c_1}),
(\ref{c_2}), (\ref{c12}):
\begin{equation}\label{dgdl1}
\lim_{\varepsilon \to+0}(G(x,\xi,l+\varepsilon)-G(x,\xi,l-\varepsilon))={1
\over[2l+1]_q}\left({1 \over c(l)}\psi_{-1-l}(\xi)+{1
\over c(-1-l)}\psi(\xi)\right)f_l(x).
\end{equation}
Now (\ref{dgdl}) follows from (\ref{dgdl1}) and lemma \ref{f_l}. The case $x
\ge \xi$ is completely similar to the case $x \le \xi$. \hfill $\Box$

\medskip

 {\sc Remark 4.15.} There is a natural generalization of proposition
\ref{dgdlst}. Let ${\rm Re}\,l=-{1 \over 2}$ and let $\gamma(\varepsilon)$
be such a parametric smooth curve on the complex plane that $\gamma(0)=l$,
${\textstyle d \gamma(0)\over \textstyle d \varepsilon}>0$. Then
\begin{equation}\label{dgdg}
\lim_{\varepsilon
\to+0}(G(x,\xi,\gamma(-\varepsilon))-G(x,\xi,\gamma(\varepsilon)))=
\frac{f_l(\xi)f_l(x)}{c(l)c(-1-l)[2l+1]_q}.
\end{equation}

\medskip \stepcounter{theorem}

 Remind that the spectrum of $\Box^{(0)}$ is simple, purely continuous and
fills a segment. This segment was parametrized as follows:
$$\lambda(l)=-\frac{(1-q^{-2l})(1-q^{2l+2})}{(1-q^2)^2},\qquad l=-{1 \over
2}+i \rho,\;0 \le \rho \le{\pi \over h}.$$
Here, as before, $h=-2 \ln \,q$. Note that
\begin{equation}\label{dldl}
{d \lambda \over dl}={1 \over(1-q^2)^2}d(q^{-2l}+q^{2l+2})={h
\over(1-q^2)^2}(q^{-2l}-q^{2l+2}).
\end{equation}

 Apply proposition \ref{specpr} to $\Box^{(0)}$. An application of
(\ref{dgdg}), (\ref{dldl}) yields the main result of this section, which was
kindly communicated to the authors by L. I. Korogodsky.

 Associate to each finitely supported function $f(x)$ on $q^{-2{\Bbb Z}_+}$
the function
$$\widehat{f}(\rho)=\int \limits_1^\infty{_3}\Phi_2
\left[{x,q^{-2l},q^{2(l+1)};q^2;q^2 \atop q^2,0}\right]f(x)d_{q^2}x$$
on the segment $\left[0,{\textstyle \pi \over \textstyle h}\right]$. Here
$l=-{1 \over 2}+i \rho$, $h=-2 \ln \,q$.

\medskip

 {\sc Example 4.16.} Let $f_0(x)=\cases{1,&$x=1$\cr 0,&$x \ne 1$}$. Then
\begin{equation}\label{f_0}
\widehat{f}_0(\rho)=1-q^2.
\end{equation}

\medskip \stepcounter{theorem}
 Remind a well known result of operator theory (\cite{AG}):

\medskip

\begin{proposition} Let $A$ be a bounded selfadjoint operator with simple
spectrum in a Hilbert space $H$, $E_t$ the spectral measure of $A$, and $g$
such a vector that the linear span of $\{A^mg \}_{m \in{\Bbb Z}}$ is dense
in $H$. With $\sigma(t)=(E_tg,g)$, the map
$$f(t) \mapsto \int \limits_{-\infty}^\infty f(t)dE_{t}g$$
is a unitary operator from $L_{\sigma}^2(-\infty,\infty)$ onto $H$.  This
unitary map sets up the equivalence of $A$ and the multiplication operator
$f(t) \mapsto tf(t)$ in $L_{\sigma}^2(-\infty,\infty)$.

\end{proposition}

\medskip

 Now one can prove the following

\begin{proposition} Consider a Borel measure \begin{equation}\label{sigma} d
\sigma(\rho)=\frac{1}{2 \pi}\cdot \frac{h}{1-q^2}\cdot \frac{d \rho}{c(-{1
\over 2}+i \rho)c(-{1 \over 2}-i \rho)} \end{equation} on the segment
$[0,{\textstyle \pi \over \textstyle h}]$. The linear operator $f \mapsto
\widehat{f}$ is extendable by a continuity up to a unitary operator $u:L^2(d
\nu)_q^{(0)} \to L^2(d \sigma)$. For all $f \in L^2(d \nu)_q^{(0)}$, $$u
\cdot \Box^{(0)}f=\lambda(l)uf.$$
\end{proposition}

\medskip

 To conclude, note that the measure $dm(l)$ could be derived from the
measure $d \sigma(\rho)$ via the substitution $l=-{1 \over 2}+i \rho$.

\bigskip

\section{Fourier transform}

 In \cite[section 5]{SSV2} a unitary operator
\begin{equation}
\overline{i}:L^2(d \nu)_q \to \bigoplus \int \limits_{{\frak
L}_0}\overline{V}^{(l)}dm(l)
\end{equation}
was constructed, with $\overline{V}^{(l)}$ being a completion of the $U_q
\frak{su}(1,1)$-module $V^{(l)}$, equipped with an invariant scalar product.
By the results of the previous section,
$$\bigoplus \int \limits_{{\frak L}_0}\overline{V}^{(l)}dm(l)\simeq
\bigoplus \int \limits_0^{\pi/h}\overline{V}^{(-{1 \over 2}+i \rho)}d
\sigma(\rho),$$
with $d \sigma$ being the measure (\ref{sigma}), and the modules $V^{(-{1
\over 2}+i \rho)}$ could be replaced by the isomorphic modules ${\Bbb
C}[z]_{q,-{1 \over 2}+i \rho}$. The linear operator $\overline{i}$ is
replaced by a completion in $L^2(d \nu)_q$ of a morphism of $U_q
\frak{sl}_2$-modules given by $if_0=1-q^2$. (This relation follows from
(\ref{f_0}); it determines unambiguously a morphism of $U_q
\frak{sl}_2$-modules by \cite[proposition 3.9]{SSV2}).

 Remind the notation (see \cite{SSV1}):
$$P_l^t=(q^2z^*\zeta;q^2)_{-l} \cdot(z \zeta^*;q^2)_{-l}(1-\zeta
\zeta^*)^l\in D(\Xi \times X)_q'.$$

\medskip

\begin{proposition} For all $f \in D(U)_q$,
$$if=\int \limits_{U_q}P_{{1 \over 2}+i \rho}^t(z,\zeta)f(\zeta)d \nu$$
\end{proposition}

\smallskip

 {\bf Proof.} It is easy to show that for all $\rho \in[0,{\textstyle \pi
\over \textstyle h}]$ the linear integral operator $i_ \rho:f \mapsto
\displaystyle \int \limits_{U_q}P_{{1 \over 2}+i \rho}^t(z,\zeta)f(\zeta)d
\nu$ maps the vector space $D(U)_q$ into ${\Bbb C}[\partial U]_{q,-{1 \over
2}+i \rho}$. Now our statement follows from the following two lemmas.

\medskip

\begin{lemma} $i_ \rho f_0=1-q^2$.
\end{lemma}

\smallskip

 {\bf Proof.} Apply the decomposition
$$P_l^t=\sum_{j>0}\zeta^j \cdot
\psi_j(\xi)+\psi_0(\xi)+\sum_{j>0}\psi_{-j}(\xi)\zeta^j,$$
described in \cite{SSV1}. It is easy to show that only the term
$\psi_0(\xi)$ contributes to the integral $i_ \rho f_0$. On the other hand,
$\psi_0(1)=1$, $\displaystyle \int \limits_{U_q}1 \cdot f_0d \nu=1-q^2$.
\hfill $\Box$

\medskip

\begin{lemma}\label{irho} The linear operator $i_ \rho:D(U)_q \to{\Bbb
C}[\partial U]_{q,-{1 \over 2}+i \rho}$ is a morphism of $U_q
\frak{sl}_2$-modules.
\end{lemma}

\smallskip

 {\bf Proof.} Consider the integral operator
$$j_ \rho:{\Bbb C}[\partial U]_{q,l}\to D(U)_q',\qquad j_ \rho:f \mapsto
\int \limits_0^{\pi/h}P_{{1 \over 2}-i \rho}(z,e^{i \theta})f(e^{i
\theta}){d \theta \over 2 \pi}.$$
It is a morphism of $U_q \frak{sl}_2$-modules, as it was noted in section 3.
Equip the $U_q \frak{su}(1,1)$-modules $D(U)_q$, ${\Bbb C}[\partial
U]_{q,l}$ with invariant scalar products
$$D(U)_q \times D(U)_q \to{\Bbb C},\qquad f_1 \times f_2 \mapsto \int
\limits_{U_q}f_2^*f_1d \nu,$$
$${\Bbb C}[\partial U]_{q,-{1 \over 2}+i \rho}\times{\Bbb C}[\partial
U]_{q,-{1 \over 2}+i \rho}\to{\Bbb C},\qquad f_1 \times f_2 \mapsto \int
\limits_{\partial U}f_2^*f_1{d \theta \over 2 \pi}.$$

 It follows from the definitions that the integral operator with a kernel
$K=\sum \limits_ik_i''\otimes k_i'$ is conjugate to the integral operator
with the kernel $K^t=\sum \limits_ik_i'^* \otimes k_i''^*$. Hence $i_
\rho=j_ \rho^*$, and $i_ \rho$ is a morphism of $U_q \frak{sl}_2$-modules
since this is the property of $j_ \rho$ (see \cite[section 5]{SSV2}). \hfill
$\Box$

\medskip

 It follows from the proof of lemma \ref{irho} that
$\overline{j}=\overline{i}^*$ is the integral operator
$$\overline{j}:f(e^{i \theta},\zeta)\mapsto \int \limits_0^{\pi/h}\int
\limits_0^{2 \pi}P_{{1 \over 2}-i \rho}(z,e^{i \theta})f(e^{i
\theta},\rho){d \theta \over 2 \pi}d \sigma(\rho).$$
Since $\overline{i}$ is unitary (see \cite{SSV2}), $\overline{j}\cdot
\overline{i}=\overline{i}\cdot \overline{j}=1$. Hence
$\overline{i},\overline{j}$ coincide with the operators $F,F^{-1}$
introduced in \cite{SSV1}, respectively. This implies the statement of
\cite[proposition 6.1]{SSV1}.

\end{document}